\documentclass[reqno]{amsart}[12pt]
\usepackage{amsmath, amssymb, latexsym, amsthm, amsfonts, mathrsfs, amscd}
\usepackage{enumerate}[1.)]
\renewcommand{\eqref}[1]{(\ref{#1})}   
\theoremstyle{plain}

\theoremstyle{remark}

\begin{document}

\title{Commutative rings with every non-maximal ideal finitely generated} 
\date{}
\author{Souvik Dey}
\address{}
\email{dey0976@gmail.com}



\maketitle

ABSTRACT : In commutative ring theory, there is a theorem of Cohen which states that if in a commutative ring all prime ideals are finitely generated then every ideal is finitely generated. However, it is known that having only maximal ideals finitely generated doesn't imply all ideals be finitely generated. In his article we ask the question that what happens if we assume all non-maximal ideals are finitely generated, and we answer the question by showing that indeed then all maximal ideals are also finitely generated i.e. the ring becomes Noetherian. 
\newline
\renewcommand{\theenumi}{(\roman{enumi})}

All the rings in this article will be commutative with unity. In the theory of commutative rings, a basic interesting and useful result of Cohen is that for a ring to be Noetherian i.e. for every ideal to be finitely generated, it is enough to have that every prime ideal be finitely generated [1]. However, it is known that only requiring finite generation of maximal ideals in a ring is not enough to force finite generation of every ideal. In this article, we consider the question that what happens if all non-maximal ideals in a ring are finitely generated. We will show that finite generation of all non-maximal ideals is indeed enough to  get finite generation of all ideals i.e. the Noetherian ness of the ring. We will extensively use the the easy basic facts that for a ring $R$, the $R$-submodules of the $R$-module $R$ are exactly the ideals and an ideal is finitely generated is equivalent to saying it is finitely generated as an $R$-module and also the fact that for an ideal $I$ of $R$, the ideals of $R/I$ are exactly the $R$-submodules of the $R$-module $R/I$. 
\newline

We first begin with the special case of integral domains, in which case we can actually prove much more:
\newline

$\mathbf{Theorem}$ 1 : Let $R$ be an integral domain. If every non-radical ideal of $R$ is finitely generated, then every ideal of $R$ is also finitely generated, i.e. the ring is Noetherian. So in particular, if every non-prime ideal in an integral domain is finitely generated, then the domain is Noetherian.  
\newline

Proof: Let $R$ be an integral domain. Let $0\ne r \in R$ and $I$ be a proper ideal of $R$. Then $r \notin rI$. Because if $r\in rI$, then $r=ri$ for some $i\in I$ and then $0\ne r$ and $R$ is a domain implies $1=i\in I$, contradicting $I$ is proper. Also, since $R$ is a domain, via the natural map sending every $j\in I$ to $rj \in rI$, we have $rI \cong I$ as $R$-modules, for every $0\ne r\in R$. Now let $0\ne J$ be any proper ideal of $R$. Let $0\ne x \in J$. Since $J$ is a proper ideal, so $x \notin xJ$. But $x^2 \in xJ$. Thus $xJ$ is not a radical ideal of $R$, hence $xJ$ is finitely generated. Then due to $xJ \cong J$ as $R$-modules, $J$ is also finitely generated. Since $J$ was any arbitrary ideal, this proves the claim. $\square$ 
\newline

Now some terminologies. Following [2], given a ring $R$, we will call an $R$-module $M$ to be almost finitely generated (a.f.g. in short) if $M$ is not finitely generated as $R$-module but every proper $R$-submodule of $M$ is finitely generated. We will use this terminology and the results of [2] and [3] freely in the proof of the following theorem which is the main theorem of this article:
\newline

$\mathbf {Theorem}$ 2 : Let $R$ be a commutative ring with unity. If every non-maximal ideal of $R$ is finitely generated, then every maximal ideal of $R$ is also finitely generated, i.e. the ring is Noetherian. 
\newline

Proof: First let us observe that we are through if $R$ is not zero dimensional i.e. if there is a non-maximal prime ideal $P$ in $R$. Indeed if $P$ is a non-maximal prime ideal, then $P$ is finitely generated, and in the integral domain $R/P$, every non-maximal ideal is finitely generated (since every non-maximal ideal of $R/P$ corresponds to a non-maximal ideal of $R$). Thus every non-prime ideal of $R/P$ is also finitely generated, and then by Theorem 1, $R/P$ is a Noetherian ring, i.e. $R/P$ is  a Noetherian $R$-module. Now $P$ is a non-maximal ideal of $R$, so every $R$-submodule of $P$ i.e. every ideal of $R$ contained in $P$ is also non-maximal,  hence finitely generated as $R$-module by hypothesis. Thus $P$ is a Noetherian $R$-module. Hence both $R/P$ and $P$ are Noetherian $R$-modules, thus $R$ is a Noetherian $R$-module i.e. $R$ is a Noetherian ring.

Now in general, let $R$ be a ring with all non-maximal ideals finitely generated and let if possible, $R$ be not Noetherian. Then by the argument above, we get that $R$ is zero-dimensional. Since $R$ is not Noetherian, by hypothesis, we get that $R$ has a maximal ideal, say $\mathfrak m$, which is not finitely generated. Now proper $R$-submodules of $\mathfrak m$ are exactly ideals of $R$ which are properly contained in $\mathfrak m$. Now any ideal properly contained in $\mathfrak m$ is not maximal, hence by hypothesis, is finitely generated as $R$-module. Thus we see that $\mathfrak m$ is an a.f.g.  $R$-module (see [2]). Hence $ann_R (\mathfrak m)$ is a prime ideal of $R$ ( [2, Proposition 1.1 (3)], also [3, Proposition 1.2 (b)] ) but is not a maximal ideal ([2, Proposition 1.3]). Thus $R$ is not zero-dimensional, contradiction ! Thus $R$ must be Noetherian.  
\newline

Now we give an alternative proof of $\mathbf {Theorem}$ 2 which avoids using any result from [2] or [3].
\newline

$\mathit{Alternative}$ $\mathit {proof}$ $\mathit {of}$ $\mathit {Theorem }$ 2: First we prove our claim for non-local rings. So let $R$ be a non-local ring whose every non-maximal ideal is finitely generated. We want to show every maximal ideal of $R$ is also finitely generated. So let $\mathfrak m$ be a maximal ideal. Since $R$ is non-local, pick a maximal ideal $\mathfrak n$ of $R$ distinct from $\mathfrak m$. Then $\mathfrak m \cap \mathfrak n $ is a non-maximal ideal , hence is finitely generated. Also $\mathfrak m + \mathfrak n=R$. Now we have the following $R$-module isomorphisms $$ \mathfrak m/ (\mathfrak m \cap \mathfrak n) \cong  (\mathfrak m + \mathfrak n)/ \mathfrak n \cong R/\mathfrak n $$  Since $\mathfrak n$ is a maximal ideal, hence $\mathfrak m/ (\mathfrak m \cap \mathfrak n)\cong R/\mathfrak n$ is a simple $R$-module, hence a cyclic $R$-module, so in particular, a finitely generated $R$-module. And   $\mathfrak m \cap \mathfrak n $ also is a finitely generated $R$-module, thus $\mathfrak m$ is a finitely generated $R$-module. Since $\mathfrak m$ was an arbitrary maximal ideal, this proves our claim. [I would like to acknowledge  the discussions in the Math overflow question as mentioned in [4] here, especially the argument by Keith Kearnes for this part of the argument]

Now for the local case, let $R$ be a local ring with unique maximal ideal $\mathfrak m$ such that every non-maximal ideal of $R$ is finitely generated. We need to show that $ \mathfrak m$ is finitely generated. If $\mathfrak m=0$, then we are done. So assume $\mathfrak m\ne 0$, hence $ann_R (\mathfrak m)$ is a proper ideal of $R$. So two cases may appear: either $ann_R (\mathfrak m) \subsetneq \mathfrak m$ or else $ann_R(\mathfrak m)= \mathfrak m$. First suppose $ann_R (\mathfrak m) \subsetneq \mathfrak m$, and fix $r \in \mathfrak m \setminus ann_R (\mathfrak m)$. Then $r \mathfrak m\ne 0$, so $ann_R (r) \neq \mathfrak m$, thus $ann_R(r)$ is a finitely generated ideal. Moreover $r\ne 0$, so $ann_R(r)$ is a proepr ideal and $(R,\mathfrak m)$ is local, hence $ann_R (r) \subsetneq \mathfrak m$. Also $r\mathfrak m\ne m$, because $r\mathfrak m=\mathfrak m$ would imply $r=rx$ for some $x\in \mathfrak m$, then $(1-x)r=0$, and since $(R,\mathfrak m)$ is local we have $1-x$ is a unit, so $r=0$, contradicting $r\notin ann_R(\mathfrak m)$. Thus $r\mathfrak m$ is non-maximal , hence finitely generated. Now we have the canonical $R$-module isomorphism (induced by the multiplication map by $r$) $\mathfrak m/ann_R(r) \cong r\mathfrak m$.  Thus $\mathfrak m /ann_R(r)$ is a finitely generated $R$-module, which along with finite generation of $ann_R(r)$ implies $\mathfrak m$ is finitely generated $R$-module. This takes care of the $ann_R (\mathfrak m) \subsetneq \mathfrak m$ case. Now for the other case, let $ann_R(\mathfrak m)=\mathfrak m$. Then $\mathfrak m$ is naturally an $R/\mathfrak m$-module, the module multiplication being given by $(x+\mathfrak m).y=xy, \forall x\in R, y\in \mathfrak m$. Thus the $R/\mathfrak m$-submodules of $\mathfrak m$ are exactly the ideal of $R$ contained in $\mathfrak m$. Now if $\mathfrak m$ were not a finitely generated ideal, then $\mathfrak m$ would be an infinite dimensional $R/\mathfrak m$-vector space. But every infinite dimensional vector space has an infinite dimensional proper subspace, hence we would have an ideal properly contained in $\mathfrak m$ which is infinite dimensional as $R/\mathfrak m$-vector space, and then due to the canonical $R/\mathfrak m$-module structure of $\mathfrak m$, such a non-maximal ideal can not be finitely generated as $R$-module, contradiction ! Thus $\mathfrak m$ must be finitely generated.  
$\square$ 
\newline

REFERENCES : 
\newline 

1. I. KAPLANSKY, Commutative rings, University of Chicago Press. 
\newline

2. W. D. WEAKLEY, Modules whose proper submodules are finitely generated, J. Algebra 84 (1983), 189 - 219.
\newline

3.  E. ARMENDARIZ, Rings with an almost Noetherian ring of fractions, Math. Scand. 41 (1977), 15 - 18. 
\newline

4. local ring all whose non-maximal ideals are finitely generated, link address: https://mathoverflow.net/questions/302468/local-ring-all-whose-non-maximal-ideals-are-finitely-generated

\end{document}